\documentstyle{article}

                    \newcommand{\w}{\widetilde}
                    \renewcommand{\l}{\lambda}
                    \newcommand{\fr}{\frac}

                    \newcommand{\bequ}{\begin{equation}}
                    \newcommand{\eequ}{\end{equation}}
                    \newcommand{\beqa}{\begin{eqnarray}}
                    \newcommand{\eeqa}{\end{eqnarray}}
                    \newcommand{\beq}{\begin{eqnarray*}}
                    \newcommand{\eeq}{\end{eqnarray*}}
										\newcommand{\barr}{\begin{array}}
\newcommand{\earr}{\end{array}}

\newcommand{\bfr}{\begin{flushright}}
\newcommand{\efr}{\end{flushright}}
\newcommand{\bfl}{\begin{flushleft}}
\newcommand{\efl}{\end{flushleft}}

\newtheorem{opr}{Definition}
\newtheorem{sle}{Corollary}
\newtheorem{teo}{Theorem}
\newtheorem{lem}{Lemma}
\newtheorem{pre}{Supposition}

\newcommand{\bo}{\begin{opr}}
\newcommand{\eo}{\end{opr}}
\newcommand{\bs}{\begin{sle}}
\newcommand{\es}{\end{sle}}
\newcommand{\bt}{\begin{teo}}
\newcommand{\et}{\end{teo}}
\newcommand{\bl}{\begin{lem}}
\newcommand{\el}{\end{lem}}
\newcommand{\bp}{\begin{pre}}
\newcommand{\ep}{\end{pre}}
\newcounter{rem}
\newcounter{pr}
\newcounter{remar}
\newcounter{cpr}

\title{A note on one inverse spectral problem.}
\author{Azamat M. Akhtyamov}

\date{}

\font\Sets=msbm10

\def\C{\hbox{\Sets C}}

\begin{document}

\maketitle

\begin{abstract}
The note contains the proof of
the uniqueness theorem for
 the inverse problem in the case
 of $n$-th order differential equation.
\end{abstract}

\vspace{0.2cm}

The inverse spectral problem is studied
in papers of many authors
(\cite{Ambarzumijan 29}--\cite{Levitan 84}).
 Extensive
bibliographies for
the inverse spectral problem
can be found in
\cite{Naymark 69 a}--
\cite{Levitan 84}.

Let's consider the spectral problem for the common differential
differential
equation:
\bequ
F(x, \, y(x),\, y'(x),\, \dots , \, y^{(n)},\, q_1(x),\,
q_2(x),\, \dots ,\, q_m(x), \l ) = 0 \label{a-s-s F} \eequ
with common boundary conditions
\bequ
U_j(y(x),\, \l ,\, a_0, \, a_1,\, \dots ,\, a_s)=0, \qquad j=1,\,\dots ,\,
n.
\label{a-s-s U}
\eequ
Here $ x\in [0,1], $ $ \l $ is eigenvalue parameter,
$ q _ i $ ($ i = 1, \, \dots \, m $) are uncertain factors of the equation,
$ a _ i $ ($ i = 0, \, \dots, \, s $) are uncertain constants of boundary
conditions, $ q _ i \in C^1[0,1]$ ($ i = 1, \, \dots \, m $),
$ a _ i \in \C $ ($ i = 0, \, \dots, \, s $).

The spectral problem defined by equalities (\ref{a-s-s F})--(\ref{a-s-s U})
we shall name $F.$

Along with the problem $ F $ we shall consider $ m $ problems:
$ A _ i $:
\beq
-y '' + q _ i (x) \, y = \l \, y, \\
y ' (0) = 0,
\\
y ' (1) = 0.
\eeq
and
one more problem $ A _ {m + 1} $:
\beqa
y '' + 3 \, y ' + 2\l ^ 2 \, y = 0, \label {a-s-s A m1} \\
y (0) = 0,
\label {a-s-s A m2} \\
y ' (1) + a (\l) \cdot y (1) = 0, \label {a-s-s A m3} \eeqa

\vspace{0.2cm}

{\bf Theorem.}

{\sl
If eigenvalues of the problems
$ A _ i $ and $ \w A _ i $ ($ i = 1, \, 2, \, \, \dots, \, m + 1 $)
coincide with their algebraic multiplicities,
then the factors of the equations and
the constant in the boundary conditions of the problems
$ F $ and $ \w F $ coincide, that is
$ q _ i(x) \equiv \w q _ i(x) ,$
$ a _ {k} = \w a _ {k} $ $ i = 1, \, 2, \, \dots, \,
m $
$ k = 1, \, 2, \, \dots, \, s. $}

\vspace{0.1cm}

{ \sl  Proof.} It follows from
Ambarzumijan's theorem (\cite{Ambarzumijan 29})
that the equality $ q _ i(x) = \w q _ i(x) $ is true.

Functions
$ y_ 1 (x, \, \l) = -e ^ {2\l x} + 2 \, e ^ {\l x}, $
$ y_ 2 (x, \, \l) = \fr {1} {\l} \, (e ^ {2\l x} - e ^ {\l x}) $  are
solutions of the differential equation (\ref {a-s-s A m1}),
satisfying
\bequ
y_ 1 (0, \, \l) = 1, \quad
y_ 1 ' (0, \, \l) = 0, \qquad
y_ 2 (0, \, \l) = 0, \quad
y_ 2 ' (0, \, \l) = 1.
\label {du m4} \eequ

Let
$ A (\l) $ be a polynomial
$ a _ {0} + a _ {1} \, \l + a _ {2} \, \l ^ 2 + \dots + a _ {s} \l ^ s. $

The eigenvalues $ \l _ i $ of the problems (\ref{a-s-s A m1})--(\ref{a-s-s
A
m3}) are the roots of a characteristic determinant, therefore
its satisfy to the following equation:
\bequ
\Delta (\l) =
\fr {1} {\l} \, \left (e ^ {2 \,\l} - e ^ {\l} \right) + a (\l) \, \left
(-e ^
{2 \,\l} + 2 \, e ^ {\l} \right) = 0.
\label {du m6}
\eequ

This function has infinite number of the radicals. (For this reason the
equation was selected by such:
$ y '' + 3 \, y ' + 2\l ^ 2 \, y = 0. $ Generally speaking it was possible
to
select
any equation having not less $s$ pairwise different nonzero eigenvalues.)

From  (\ref {du m4}) we have
\bequ
1\cdot a _ {0} + \l _ i\cdot a _ {1} + \l _ i ^ {\, 2} \cdot a _ {2} +
\dots +
\l _ i ^ {\, s} \cdot a _ {s}
= - \fr {e ^ {2 \,\l _ i} - e ^ {\l _ i}} {-\l _ i \, e ^ {2 \,\l _ i} + 2
\,\l
_ i \, e ^ {\l _ i}}. \label {du m7} \eequ
The equalities (\ref {du m7}) is a system of $ (s + 1) $ linear equaions
having
$ (s + 1) $
unknown
$ a _ {0}, $ $ a _ {1}, $ $ a _ {2}, $ $ \dots, $ $ a _ {s} $ of the
 boundary condition.
The determinant of this system is the Vandermonde determinant. As the
eigenvalues $
\l _ i $ are pairwise  different and are not equal to zero, the Vandermonde
determinant is not equal to zero.
Therefore system (\ref {du m7}) has a unique solution.
From here
Follows, that the constants $ a _ {0}, $ $ a _ {1}, $ $ a _ {2}, $ $ \dots,
$ $
a _ {s}$ are determined univalently.

The theorem is proved.

\end{document}